 \documentclass[11pt,a4paper]  
 {article}
\usepackage{amssymb,latexsym, amsmath, bm, amsfonts, tikz,fancyhdr,sectsty}
\usepackage{pgfplots}
\usetikzlibrary{patterns}

\usepackage{lettrine}
\newtheorem{theorem}{Theorem}[section]
\newtheorem{lemma}[theorem]{Lemma}

 \title{\bf Generalised differences and     a class of multiplier 
   operators  in Fourier analysis}
\author{Rodney Nillsen \\
\small{School of Mathematics and Applied Statistics}\\
\small{University of Wollongong,
  Northfields Avenue, Wollongong}\\
\small{New South Wales, 2522, AUSTRALIA}\\
\small{email address: nillsen@uow.edu.au}\\}
 \date{\ }
  
\begin{document}
\maketitle
 \begin{abstract}
Any zeros in the multiplier of an operator impose a  condition on a  function for it to be  in the range of the operator.  But if  each  function in a certain family $\mathcal F$ of functions  satisfies such a condition, when is  $\mathcal F$ the range of the operator?  Let $\alpha,\beta\in {\mathbb Z}$, for  $g\in L^2([0,2\pi])$ let ${\widehat g}$ be the sequence of Fourier coefficients  of $g$, and let $D$ denote differentiation. We consider the operator  $D^2-i(\alpha+\beta)D-\alpha\beta$ on the second order Sobolev space of $L^2([0,2\pi])$. The multiplier of this operator is  $-(n-\alpha)(n-\beta)$, so that ${\widehat g}(\alpha)={\widehat g}(\beta)=0$  for any function $g$ in the range of the operator.  Let $\delta_x$ denote the Dirac measure at $x$, and  let $\ast$ denote convolution. If $b\in [0,2\pi]$  let $\lambda_b$ be  the measure   
\[\footnotesize{\frac{1}{2}\left[e^{ib\left(\frac{\alpha-\beta}{2}\right)}+e^{-ib\left(\frac{\alpha-\beta}{2}\right) } \right]\delta_0- \frac{1}{2}\left[\  e^{ib\left(\frac{\alpha+\beta}{2}\right)}\,\delta_{ b}+e^{-ib\left(\frac{\alpha+\beta}{2}\right)}\,\delta_{- b}\right]}.\]
 A function of the form $\lambda_b \ast f$ is called a \emph{generalised difference}, and we let 
  ${\mathcal F}$ be the family of functions $h$ such  that $h$ is some finite sum of generalised differences.      It is shown   that   ${\mathcal F}$   is a closed subspace of $L^2({\mathbb T})$ that equals the range of  $D^2-i(\alpha+\beta)D-\alpha\beta$, and that every function in ${\mathcal F}$ is a sum of five  generalised differences. The methods use partitions of $[0,\pi/2]$  and estimates of   integrals in Euclidean space. There are applications to  the automatic continuity of linear forms.  
  \end{abstract}
\let\thefootnote\relax\footnote{2010 \emph{Mathematics Subject Classification}.  Primary 42A16, 42A45, \hfill\break Secondary 43A15, 46H40

\emph{Key words and phrases.}  Circle group, Fourier analysis, $L^2$ spaces, multiplier \hfill\break operators, generalised differences,   compact abelian groups, automatic continuity}

 \section{Introduction}
 \setcounter{equation}{0}
 The circle group  $\{z: |z|=1\}$ is denoted by ${\mathbb T}$. The mapping $x\longmapsto e^{ix}$ from $[0,2\pi)$  to ${\mathbb T}$ means that we can and will identify  ${\mathbb T}$ with $[0,2\pi)$ or $[0,2\pi]$ in the usual way, and both settings will be used (see the comments in \cite[page 1034]{ross1}).  The group operation $+$ on $[0,2\pi)$ is written additively and is  the usual addition if $0\le x+y<2\pi$, and is $x+y-2\pi$ if $2\pi \le x+y$. The space $L^2({\mathbb T})$ is the Hilbert space of square integrable complex functions on $\mathbb T$, and  $M({\mathbb T})$ denotes the set of bounded, regular, complex Borel measures on ${\mathbb T}$.   Let $\mathbb Z$ denote the set of integers, and let $\mathbb N$ denote the set of positive integers. The convolution operation in $M({\mathbb T})$ is denoted by $\ast$. Thus, if $\mu\in M({\mathbb T})$ and $n\in {\mathbb N}$, $\mu^n$ denotes $\mu\ast\mu\ast\cdots\ast\mu$, where $\mu$ appears $n$ times.   Considering   $n\in {\mathbb Z}$, $f\in L^2({\mathbb T})$ and   $\mu\in M({\mathbb T})$, we define the   Fourier coefficients ${\widehat f}(n)$  and ${\widehat {\mu}}(n)$ by
\[{\widehat f}(n)=\frac{1}{2\pi}\int_0^{2\pi}f(t)\,e^{-int}\,dt \ {\rm and}\ {\widehat {\mu}}(n)=\int_0^{2\pi}e^{-int}\,d\mu(t).\]
Note that for $\mu,\nu\in M({\mathbb T})$, $( {\mu\ast \mu})\,{\widehat {\ }}={\widehat \mu}\,{\widehat \nu}$.
Letting $\delta_x$ denote the Dirac measure at $x$, we see that ${\widehat {\delta_x}}(n)=e^{-inx}$ for $x\in [0,2\pi]$, and ${\widehat {\delta_z}}(n)=z^{-n}$ for $z\in {\mathbb T}$. For $s=0,1,2,\ldots $ the Sobolev space $W^s({\mathbb T})$ is defined by 
\[W^s({\mathbb T})=\Bigl\{f:f\in L^2({\mathbb T}) \ \, {\rm and}\, \sum_{n=-\infty}^{\infty}|n|^{2s}|{\widehat f}(n)|^2<\infty\Bigr\}.\]
The differential operator $D$  maps $W^1({\mathbb T})$ into $L^2({\mathbb T})$ and  can be defined by the property
\[D(f){\widehat {\ }}(n)=i n{\widehat f}(n), \ {\rm for \ all}\  n\in {\mathbb Z}.\]
We say that $D$ is a \emph{multiplier operator on} $W^1({\mathbb T})$ with \emph{multiplier} $in$.
Note that this means that if $f\in W^1({\mathbb T})$, $D(f){\widehat {\ }}(0)=0$. Conversely, it  is easy to see that    if $g\in L^2({\mathbb T})$ and ${\widehat g}(0)=0$, then $g=D(f)$ for some $f\in W^1({\mathbb T})$. Similarly, $D^s$ is a multiplier operator from $W^s({\mathbb T})$ into $L^2({\mathbb T})$  with multiplier $(in)^s$. 
Let $\alpha,\beta \in {\mathbb Z}$ and $s\in {\mathbb N}$ be given, and let $I$ denote the identity operator.  In this paper, and for $s\in {\mathbb N}$, we are concerned with operators  $(D^2-i(\alpha+\beta)D-\alpha\beta I)^s$ that  map $W^{2s}({\mathbb T})$ into $L^2({\mathbb T})$. These are multiplier operators with multipliers $ (-1)^s(n-\alpha)^s(n-\beta)^s$. That is, for all $f\in W^s({\mathbb T})$,
\[\bigl[(D^2-i(\alpha+\beta)D-\alpha\beta I)^s(f)\bigr]\,{\widehat {\ }}\,(n)=(-1)^s(n-\alpha)^s(n-\beta)^s{\widehat f}(n),\]
for all $n\in {\mathbb Z}$.
Note that  the zeros of the multiplier of $(D^2-i(\alpha+\beta)D-\alpha\beta I)^s$ are $\alpha$ and $\beta$.  The operator $(D^2-i(\alpha+\beta)D-\alpha\beta I)^s$ acting on $f\in W^{2s}({\mathbb T})$ eliminates the frequencies $\alpha$ and $\beta$ from $f$.

 A particular case of the above is when $\alpha=-\beta$. Then the operator $(D^2-i(\alpha+\beta)D-\alpha\beta I)^s$ becomes $(D^2+\alpha^2I)^s$, and it has the multiplier  $-n^2+\alpha^2$, which has zeros at $-\alpha$ and $\alpha$. In order to give a feeling for the ideas in this paper, here is a statement of a special case of what is perhaps  the main result.
 \emph{Let $\alpha\in {\mathbb Z}.$  Then the  following conditions (i), (ii) and (iii) are equivalent  for a function $f\in L^2({\mathbb T})$}. 
 
 \emph{(i)} ${\widehat f}(-\alpha)={\widehat f}(\alpha)=0$.   
 
 \emph{(ii)}   \emph{$f$ is a sum of five functions,  for each one $h$ of which there are  $z\in {\mathbb T}$ and  $g\in L^2({\mathbb T})$ such that  }
 \begin{equation}
 h=(z^{\alpha}+z^{-\alpha})g -(\delta_{z}+\delta_{z^{-1}})\ast g.\label{eq:differences1}
 \end{equation} 
 
\emph{ (iii) There is $g\in W^2({\mathbb T})$ such that $(D^2+\alpha^2I)(g)=f$.}

Now a function of the form $g-\delta_z\ast g$ might be called a first order difference, and one of the form $(\delta_1-\delta_z)^s\ast g$ might be called a  difference of order $s$ (see \cite{nillsen1,nillsen2}). So, a function as in (\ref{eq:differences1})  we call a type of \emph{generalised difference}.  The form of generalised differences  is derived from the  factors  of the multiplier of the operator, rather than by replacing the elements  $D$ and $D^2$ in the operator by first and second order differences.

The present ideas in proving a result of the above type involve looking at  the structure of  partitions of $[0,\pi/2]$ associated with the zeros in $[0,\pi/2]$ of the functions $\sin((n-\alpha)x)$ and $\sin((n-\beta)x)$, and estimating integrals in ${\mathbb R}^{m}$ over  sets that are Cartesian products of sets in a refined   partition.

The results here are related to work of Meisters and Schmidt \cite{meisters1}, where they proved the following result.
\emph{The  following conditions (i) and (ii) are equivalent for a function $f\in L^2({\mathbb T})$}.
 
  \emph{(i)  ${\widehat f}(0)=0$.}

\emph{(ii) $f$ is a sum of three functions, for each one $h$ of which there are  $z\in {\mathbb T}$ and  $g\in L^2({\mathbb T})$ such that}
$h=g -\delta_z\ast g.$

They deduced from this that if $T$ is a linear form on $L^2({\mathbb T})$ such that $T(f)=T(\delta_z\ast f)$ for all $z\in {\mathbb T}$ and $f\in L^2({\mathbb T})$, then $T$ is continuous.  That is,  any translation invariant linear form on $L^2({\mathbb T})$ is automatically continuous, and this was proved more generally  for compact, connected abelian groups.   Further work relating to the original results of Meisters and Schmidt  can be found, for example,  in \cite{bourgain1,  johnson1, meisters2,  meisters3,  nillsen1, nillsen2}  where there are further references. One way to think of the ideas in \cite{meisters1} is that they are concerned with the range of the differential operator $D$ whose multiplier is a linear polynomial. On the other hand, the present work is concerned with differential operators whose multipliers are   quadratic in nature. The work here goes back to \cite{meisters1}, but takes some of the ideas there in a different direction from the mainstream of later work. There are  also applications to automatic continuity of linear forms. A suitable general reference on classical Fourier series is \cite{edwards1},  and for abstract harmonic analysis on groups see \cite{hewitt1, ross1}.

\section{Background and statement of the main result}
 \setcounter{equation}{0}
 There will be cause to consider  series    of the form $\sum_{j=-\infty}^{\infty}a_j/b_j$, where $a_j,b_j\ge 0$ for all $j$. In such a case, we write $\sum_{j=-\infty}^{\infty}a_j/b_j=\infty$ if there is a term of the form $a_j/0$ with $a_j>0$. If there are any terms of the form $0/0$ we either neglect them in the sum or make them equal to $0$.  We write respectively $\sum_{j=-\infty}^{\infty}a_j/b_j<\infty$ or  $\sum_{j=-\infty}^{\infty}a_j/b_j= \infty$ if the series $\sum_{j=-\infty, a_j>0}^{\infty}a_j/b_j$ converges or diverges in the usual sense.   The following  result is due  essentially to Meisters  and Schmidt \cite[page 413]{meisters1}. 
 \begin{theorem} \label{theorem:characterisation} Let $f\in L^2([0,2\pi])$ and let $\mu_1,\mu_2,\ldots,\mu_r\in M([0,2\pi])$. Then the following conditions (i) and (ii) are equivalent.
 
 (i) There are   $f_1,f_2,\ldots, f_r\in L^2([0,2\pi])$ such that 
 $f=\sum_{j=1}^r\mu_j\ast f_j.$
 \vskip0.2cm
(ii) \hskip 3.6cm$\displaystyle \sum_{n=-\infty}^{\infty}\,\frac{|{\widehat f }(n)|^2}{\displaystyle{\sum_{j=1}^r}|{\widehat u_j}(n)|^2}\ <\ \infty.$
\end{theorem}

{\bf Proof.} This is essentially proved in \cite[pages 411-412 ]{meisters1}, but see also \cite[pages 77-88]{nillsen1} and \cite[page 23]{nillsen2}.  A more accessible proof  for the present context is on the world wide web \cite{nillsen3}.\hfill$\square$

In \cite{meisters1}, the measures  in Theorem \ref{theorem:characterisation} were taken to be of the form $\mu_b=\delta_0-\delta_b$, in which case $\mu_b([0,2\pi])=0$ and ${\widehat {\mu_b}}(n)=1-e^{-ibn}$.
In the context here, we let $\alpha,\beta\in {\mathbb Z}$ and we  will apply Theorem \ref{theorem:characterisation} using measures $\lambda_b$, $b\in [0,2\pi)$, where
\begin{equation}\lambda_b=\frac{1}{2}\left[e^{ib\left(\frac{\alpha-\beta}{2}\right)}+e^{-ib\left(\frac{\alpha-\beta}{2}\right) } \right]\delta_0- \frac{1}{2}\left[\  e^{ib\left(\frac{\alpha+\beta}{2}\right)}\,\delta_{ b}+e^{-ib\left(\frac{\alpha+\beta}{2}\right)}\,\delta_{- b}\right].\label{eq:lambdasubb}
\end{equation}
Note the the Fourier transform ${\widehat {\lambda_b}}$ of $\lambda_b$ is given for $n\in {\mathbb Z}$ by
\begin{equation}
{\widehat {\lambda_b}}(n)=\cos \left(\left(\frac{\alpha-\beta}{2}\right) b\right)-\cos \left(\left(n-\frac{\alpha+\beta}{2}\,\right)b\right). \label{eq:Fouriertransform}
\end{equation}
Consequently,
\begin{equation} {\widehat {\lambda_b}}(\alpha)={\widehat {\lambda_b}}(\beta)=0, \ {\rm for \ all}\ b\in [0,2\pi].\label{eq:fourierzero}
\end{equation}
So if $b\in [0,2\pi]$ and $g\in L^2([0,2\pi])$,  
\begin{equation}
\lambda_b\ast g=\frac{1}{2}\left[e^{ib\left(\frac{\alpha-\beta}{2}\right)} +e^{-ib\left(\frac{\alpha-\beta}{2}\right)}\right]g- \frac{1}{2}\left[\  e^{ib\left(\frac{\alpha+\beta}{2}\right)}\,\delta_{ b}+e^{-ib\left(\frac{\alpha+\beta}{2}\right)}\,\delta_{- b}\right]\ast g,\label{eq:generaliseddifference}
\end{equation}
and  we see that if $f\in L^2([0,2\pi])$ is a function of the form $\lambda_b\ast g$, then ${\widehat f}(\alpha)={\widehat f}(\beta)=0$.  A function of the form $\lambda_b\ast g$, as in (\ref{eq:generaliseddifference}) above, is called a \emph{generalised difference} in $L^2([0,2\pi])$.

The following  is a crucial technical result. It is not proved here, but is a consequence of  results in  subsequent sections.  
 
\begin{lemma} \label{lemma:cosine estimate} Let $\alpha,\beta\in {\mathbb Z}$ and let $s\in {\mathbb N}$ be given.  Then there is $M>0$ such that  for all $n\in {\mathbb Z}$ with $n\notin \{\alpha,\beta\}$,
\[\int_{{[0,2\pi]}^{4s+1}}\frac{dx_1dx_2\cdots dx_m}{ \displaystyle \sum_{j=1}^{4s+1}\left|\cos\left(\left(\frac{\alpha-\beta}{2}\right)x_j\right)-\cos \left(\left(n-\frac{\alpha+\beta}{2}\right)x_j\right)\right|^{2s}}\,\le \ M.\]

\end{lemma}

 The  central result proved here using Lemma \ref{lemma:cosine estimate}  is the following.
\begin{theorem}\label{theorem:main2} Let $\alpha,\beta\in {\mathbb Z}$.     Then the following conditions (i), (ii), (iii)  and (iv) on a function $f\in L^2([0,2\pi])$ are equivalent.  

(i)  ${\widehat f}(\alpha)={\widehat f}(\beta)=0$.

(ii) There are $m,s\in {\mathbb N}$,  $b_1,b_2,\ldots, b_m\in [0,2\pi]$ and $f_1,f_2,\ldots,f_m\in  L^2([0,2\pi])$ such that $f$ is equal to 
\begin{equation}
 \sum_{j=1}^m
 \left[
 \left(e^{ib_j\left(\frac{\alpha-\beta}{2}\right)}+e^{-ib_j\left(\frac{\alpha-\beta}{2}
 \right)}\right)
 \delta_0\,
  -\left(e^{ib_j\left(\frac{\alpha+\beta}{2}\right)}\,\delta_{b_j} 
   +e^{-ib_j\left(\frac{\alpha+\beta}{2}\right)}\,\delta_{- b_j}\right)
   \right]^s
   \ast f_j.\label{eq:differencesum}
   \end{equation}

(iii) There are $s\in {\mathbb N}$,   $b_1,b_2,\ldots, b_{4s+1}\in [0,2\pi]$ and $f_1,f_2,\ldots,f_{4s+1}\in  L^2([0,2\pi])$ such that $f$ is equal to
   \begin{equation}
 \sum_{j=1}^{4s+1}
 \left[
 \left(e^{ib_j\left(\frac{\alpha-\beta}{2}\right)}+e^{-ib_j\left(\frac{\alpha-\beta}{2}
 \right)}\right)
 \delta_0\,
  -\left(e^{ib_j\left(\frac{\alpha+\beta}{2}\right)}\,\delta_{b_j} 
   +e^{-ib_j\left(\frac{\alpha+\beta}{2}\right)}\,\delta_{- b_j}\right)
   \right]^s
   \ast f_j. \label{eq:differencesumb}
   \end{equation}

  (iv)  There are $s\in {\mathbb N}$ and $g\in W^{2s}([0,2\pi])$ such that 
  \begin{equation}\big(D^2-i(\alpha+\beta)D-\alpha\beta I\big)(g)=f.\label{eq:differentialoperator}\end{equation}
  
 \noindent When  the equivalent  conditions (i), (ii), (iii) and (iv) are satisfied and $s\in {\mathbb N}$ is given,   we have that for almost all $(b_1,b_2,\ldots, b_{4s+1})\in [0,2\pi]^{2s+1}$, there are  $f_1,f_2,\ldots,f_{4s+1}\in  L^2([0,2\pi])$ such that (\ref{eq:differencesumb}) holds. Also, the functions in $L^2([0,2\pi])$  that can be written in the form (\ref{eq:differencesumb}) form a closed vector subspace of $L^2([0,2\pi])$.
\end{theorem}
{\bf Proof.} It is obvious that (iii) implies (ii). Now, (ii) implies (i) since, as noted in  (\ref{eq:fourierzero}), the Fourier coefficients of any function appearing within the sum in (\ref{eq:differencesum})  vanish at $\alpha$ and $\beta$. 

In order to  prove that (i) implies (iii), let ${\widehat f}(\alpha)={\widehat f}(\beta)=0$, and let $M$ be the constant as in Lemma \ref{lemma:cosine estimate}.  If we integrate the function that maps $ (x_1,x_2,\ldots, x_{4s+1})\in{\mathbb R}^{4s+1}$ into 
\begin{equation} \sum_{{n=-\infty}\atop{n\ne \alpha,\beta}}\frac{ |{\widehat f}(n)|^2 }{ \displaystyle\sum_{j=1}^{4s+1}\left|\cos\left(\left(\frac{\alpha-\beta}{2}\right)x_j\right)-\cos \left(\left(n-\frac{\alpha+\beta}{2}\right)x_j\right)\right|^{2s}}\label{eq:sumintegral}
\end{equation}
over $[0,2\pi]^{4s+1}$, and interchange the order of integration and summation we obtain
\begin{align*}
&\sum_{{n=-\infty}\atop{n\ne \alpha,\beta}}^{\infty}\left(\int_{[0,2\pi]^{4s+1}}\frac{dx_1dx_2\cdots dx_{4s+1}} { \displaystyle\sum_{j=1}^{4s+1}\left|\cos\left(\left(\frac{\alpha-\beta}{2}\right)x_j\right)-\cos \left(\left(n-\frac{\alpha+\beta}{2}\right)x_j\right)\right|^{2s}}\right)  |{\widehat f}(n)|^2 \\
&\le M\sum_{n=-\infty}^{\infty}|{\widehat f}(n)|^2\\
&<\infty.
\end{align*}
We deduce that for almost all $(x_1,x_2,\ldots,x_{4s+1})\in [0,2\pi]^{4s+1}$, the sum in (\ref{eq:sumintegral}) is finite.  Using (\ref{eq:Fouriertransform}), we see that for almost all 
$(x_1,\ldots,x_{4s+1}) \in [0,2\pi]^{4s+1}$,
\begin{align*}
&\sum_{n=-\infty}^{\infty}\frac{|{\widehat f}(n)|^2}{\displaystyle\sum_{j=1}^{4s+1}|{\widehat {\lambda_{x_j}^s}}  (n)|^2} \hskip7cm
\end{align*}
\begin{align*}
&=\sum_{n=-\infty}^{\infty}\frac{{\widehat f}(n)|^2}{\displaystyle\sum_{j=1}^{4s+1}|{\widehat \lambda}_{x_j}(n)|^{2s}} \\
& =\sum_{n=-\infty}^{\infty}\frac{{\widehat f}(n)|^2}{\displaystyle\sum_{j=1}^{4s+1}|{\widehat \lambda}_{x_j}(n)|^{2s}}\\
&=\sum_{{n=-\infty}\atop{n\ne \alpha,\beta}}^{\infty}\,\frac{ |{\widehat f}(n)|^2 }{\displaystyle \sum_{j=1}^{4s+1}\left|\cos\left(\left(\frac{\alpha-\beta}{2}\right)x_j\right)-\cos \left(\left(n-\frac{\alpha+\beta}{2}\right)x_j\right)\right|^{2s}}\\
&<\,\infty.
\end{align*}
It now follows from (\ref{eq:Fouriertransform}), and the equivalence of (i) and (ii) in Theorem \ref{theorem:characterisation},  that 
 for almost all $(b_1,b_2,\ldots b_{4s+1})\in [0,2\pi]^{4s+1}$ there are $f_1,f_2,\ldots,f_{4s+1}\in  L^2([0,2\pi])$ such that 
 \[f=\sum_{j=1}^{4s+1}\lambda_{b_j}^s\ast f_j. \]
  We see, using (\ref{eq:lambdasubb}), that (\ref{eq:differencesumb}) and hence (iii) hold.
  
  Now to see that (iv) implies (i), we observe that the multiplier of \hfill\break $(D^2-i(\alpha+\beta)D-\alpha\beta I)^s$ is $(-1)^s(n-\alpha)^s(n-\beta)^s$.  Thus, if \hfill\break $f=(D^2-i(\alpha+\beta)D-\alpha\beta I)^s(g)$, as in (\ref{eq:differentialoperator}),  ${\widehat f}(\alpha)={\widehat f}(\beta)=0$.
  
  Now, we prove (i) implies (iv). If  ${\widehat f}(\alpha)={\widehat f}(\beta)=0$, define $g$ in terms of ${\widehat g}$ by putting ${\widehat g}(\alpha)={\widehat g}(\beta)=0$ and ${\widehat g}(n)=(-1)^s{\widehat f}(n)/(n-\alpha)^s(n-\beta)^s$, for $n\ne \alpha$ and $n\ne \beta$. Then $\sum_{n=-\infty}^{\infty}|n|^{4s}|{\widehat g}(n)|^2<\infty$, so $g\in W^{2s}([0,2\pi])$. Also, as the multiplier of $(D^2-i(\alpha+\beta)D-\alpha\beta I)^s$ is $(-1)^s(n-\alpha)^s(n-\beta)^s$, we see that $(D^2-i(\alpha+\beta)D-\alpha\beta I)^s(g)\,{\widehat {\ }}={\widehat f}$ and so  $(D^2-i(\alpha+\beta)D-\alpha\beta I)^s(g)=f$.
  
  Finally, that the functions expressible as in  (\ref{eq:differencesumb}) form a closed subspace follows from the fact that such functions are characterised by (i).
  ${\ }$\hfill$\square$

 \section{ Partitioning, zeros and inequalities}
\setcounter{equation}{0}

In this section we make observations and  obtain results for the later proving of Lemma \ref{lemma:cosine estimate}. 
 \begin{lemma} \label{lemma:integralsequal} Let $m,s\in {\mathbb N}$ and  $\alpha,\beta,
n\in {\mathbb Z}$.  Then, 
\begin{align}&\int_{{[0,2\pi]}^m}\frac{dx_1dx_2\cdots dx_m}{\displaystyle\sum_{j=1}^{m}\left|\cos\left(\left(\frac{\alpha-\beta}{2}\right)x_j\right)-\cos \left(\left(n-\frac{\alpha+\beta}{2}\right)x_j\right)\right|^{2s}}\nonumber\\
&=2^{2m-2s}\int_{{[0,\pi/2]}^m}\frac{dx_1dx_2\cdots dx_m}{\displaystyle\sum_{j=1}^m\sin^{2s}((n-\alpha)  x_j)\sin^{2s}((n-\beta) x_j)}, \label{eq:int1}
\end{align}
but note that both integrals may be infinite.
\end{lemma}

{\bf Proof.}
The use of a familiar trigonometric formula and a simple substitution gives
\begin{align}
&\int_{{[0,2\pi]}^m}\frac{dx_1dx_2\cdots dx_m}{\displaystyle\sum_{j=1}^{m}\left|\cos\left(\left(\frac{\alpha-\beta}{2}\right)x_j\right)-\cos \left(\left(n-\frac{\alpha+\beta}{2}\right)x_j\right)\right|^{2s}}\nonumber\\
& =2^{m-2s}\int_{[0,\pi]^m}\frac{dx_1dx_2\cdots dx_m}{\displaystyle\sum_{j=1}^m\sin^{2s}((n-\alpha) x_j)\sin^{2s}((n-\beta) x_j)}.\label{eq:int2}
\end{align}
Note that for $x\in {\mathbb R}$ and $\ell\in {\mathbb Z}$ , $|\sin(\ell x)|=|\sin (\ell(\pi-x))|$. Also, note that $[0,\pi)^m$ is the disjoint union of the $2^m$ sets $\prod_{t=1}^mJ_t$ where, for each $t$, $J_t$ is either $[0,\pi/2)$ or $ [\pi/2,\pi)$. Using  substitutions for $\pi-x_j$, as needed in (\ref{eq:int2}), we see that  (\ref{eq:int1}) follows from (\ref{eq:int2}).
${\ }$ \hfill$\square$

We are aiming to estimate the integral in Lemma \ref{lemma:integralsequal}. Motivated by (\ref{eq:int1}), we consider the zeros in $[0,\pi/2]$ of $\sin ((n-\alpha)x)$ and $\sin((n-\beta)x)$. Some preliminary notions are needed. 

{\bf Definitions.} If $J$ is an interval we denote its length by $\lambda(J)$. Let $[a,b]$ be a closed interval with $\lambda([a,b])>0$. A family $\{J_0,J_1,\ldots,J_{r-1}\}$ of closed intervals having non-empty interiors  is a \emph{partition}  of $[a,b]$ if $\cup_{j=0}^{r-1}J_j=[a,b]$ and any two intervals in the family  have at most a single point in common. In such a case, the intervals may be arranged so that the  right endpoint of $J_{j-1}$ is the left endpoint of $J_j$ for all  $j=1,2, \ldots,r-1$.  Note that in the sense used here, the sets in a partition  are not pairwise disjoint.  \begin{lemma}\label{lemma:partitionsa} Let $[a,b]$ be a closed interval with $\lambda([a,b])>0$. Let $R_0,R_1,\ldots,$\hfill $R_{r-1}$ be  closed intervals in a partition ${\mathcal P}_1$ of $[a,b]$. Let $S_0,S_1,\ldots,S_{s-1}$  be   closed intervals in a partition ${\mathcal P}_2$ of $[a,b]$. Let 
\[{\mathcal A}=\big\{(j,k):  0\le j\le r-1, 0\le k\le s-1\ \hbox{and}\ \lambda(R_j\cap S_k)>0 \bigr\},\]
and put 
\begin{equation}{\mathcal P}=\bigl\{R_j\cap S_k: (j,k)\in {\mathcal A} \bigr\}.\label{eq:refinement}
\end{equation}
Then, ${\mathcal P}$ is a partition of $[a,b]$ and  the number of intervals in   ${\mathcal P}$ is at most 
\begin{equation}r+s-1.\label{eq:partitionestimate1}
\end{equation}
\end{lemma}

Proof. It is clear that ${\mathcal P}=\{R_j\cap S_k: (j,k)\in {\mathcal A}\}$  as in (\ref{eq:refinement}) is a partition of $[a,b]$. Given any two partitions ${\mathcal P}_1$, ${\mathcal P}_2$ as in the Lemma, we will denote the partition given as in (\ref{eq:refinement}) by ${\mathcal P}({\mathcal P}_1,{\mathcal P}_2)$.

We proceed by induction on $r$. If $r=1$, ${\mathcal P}_1=\{R_0\}=\{[a,b]\}$, so that ${\mathcal P}({\mathcal P}_1,{\mathcal P}_2)={\mathcal P}_2$ and ${\mathcal P}({\mathcal P}_1,{\mathcal P}_2)$ has $s$ elements. In this case, $s=r+s-1$, and we see that when $r=1$, (\ref{eq:partitionestimate1}) holds for all $s\in {\mathbb N}$.  Similarly, it is easy to see that when $r=2$,  (\ref{eq:partitionestimate1}) holds for all $s\in {\mathbb N}$.

So, let $r\in {\mathbb N}$ with $r\ge 3$ be such that (\ref{eq:partitionestimate1}) holds for all partitions ${\mathcal P}_1$ having $r$ intervals and for all partitions ${\mathcal P}_2$ having any number of intervals. Let ${\mathcal P}_3=\{Q_0, Q_1,\ldots,Q_r\}$ be a partition of $[a,b]$ into $r+1$ intervals and let ${\mathcal P}_2=\{S_0,S_1,\ldots,S_{s-1}\}$ be a partition of $[a,b]$ into $s$ intervals. Consider the partition ${\mathcal P}_1=\{Q_0\cup Q_1,Q_2\ldots, Q_r\}$, which has $r$ intervals.  By the inductive hypothesis,  ${\mathcal P}({\mathcal P}_1,{\mathcal P}_2)$ has at most $r+s-1$ intervals. Let $\xi$ be the right endpoint of $Q_0$, which is also the left endpoint of $Q_1$.  Now, in passing from ${\mathcal P}_1$ to ${\mathcal P}_3$, we do so by using $\xi$ to divide the single interval $Q_0\cup Q_1$ into the two intervals $Q_0$ and $Q_1$. If $\xi$ is not an endpoint of any interval in ${\mathcal P}_2$, this will divide one  interval in ${\mathcal P}({\mathcal P}_1,{\mathcal P}_2)$ into two subintervals belonging to ${\mathcal P}({\mathcal P}_2,{\mathcal P}_3)$. If $\xi$ is an endpoint of some interval in ${\mathcal P}_2$, this division does not  increase the number of intervals  in ${\mathcal P}({\mathcal P}_2,{\mathcal P}_3)$ in going from    ${\mathcal P}_1$ to  ${\mathcal P}_3$. In either case, we see that ${\mathcal P}({\mathcal P}_2,{\mathcal P}_3)$ has at most 
\[r+s-1+1=r+s=(r+1)+s-1\]
intervals, showing that (\ref{eq:partitionestimate1}) holds with $r+1$ in place of $r$. The result follows by induction.\hfill$\square$

  {\bf Definition.} Let  $[a,b]$ be a closed interval of positive length. Let ${\mathcal P}_1=\{R_0,R_1,\ldots,R_{r-1}\}$  and ${\mathcal P}_2=\{S_0,S_1,\ldots,S_{s-1}\}$  be   two partitions of $[a,b]$.   Let ${\mathcal P}$ be the partition of $[a,b]$ as given by (\ref{eq:refinement}) in Lemma  \ref{lemma:partitionsa}.
  Then   ${\mathcal P}$  is called the \emph{refinement} of the partitions ${\mathcal P}_1$ and ${\mathcal P}_2$.

Now, let $n,\gamma\in {\mathbb Z}$ with $n\ne \gamma$ be given. We construct  an associated  partition ${\mathcal P}(\gamma)$ of $[0,\pi/2]$,  as follows.

(i)  When $|n-\gamma|$ is even we define $(|n-\gamma|+2)/2$  closed subintervals 
$Q_0,Q_1,\ldots, Q_{|n-\gamma|/2}$ of $[0,\pi/2]$ by putting 
\begin{equation}Q_0=\left[0,\frac{\pi}{2|n-\gamma|}\right ], Q_{|n-\gamma|/2}=\left[\frac{\pi(|n-\gamma|-1)}{2 |n-\gamma|} ,\frac{\pi}{2}\right],\ {\rm and}\nonumber\end{equation}
\begin{equation}\ Q_j=\left[\frac{\pi(j-1/2)}{|n-\gamma|},\frac{\pi(j+1/2)}{|n-\gamma|}\right], \label{eq:Rjdefinition1}
\end{equation}
for $j=1,2 \ldots, (|n-\gamma|-2)/2$.

(ii)  When $|n-\gamma|$ is odd we define $(|n-\gamma|+1)/2$  closed subintervals 
$Q_0,Q_1,\ldots, Q_{(|n-\gamma|-1)/2}$ of $[0,\pi/2]$ by putting     \begin{equation}Q _0=\left[0,\frac{\pi}{2|n-\gamma|}\right ]\ {\rm and}\ Q_j=\left[\frac{\pi(j-1/2)}{|n-\gamma|},\frac{\pi(j+1/2)}{|n-\gamma|}\right], \label{eq:Rjdefinition2}
\end{equation}
the latter for $j=1,2 \ldots, (|n-\gamma|-1)/2$.

Put $\theta(\hbox{$r$})=(r+2)/2$ if $r$ is even,  and $\theta(\hbox{$r$})=(r+1)/2$ if $r$ is odd. Then, with $n$ and $\gamma$ as given, (\ref{eq:Rjdefinition1}) and (\ref{eq:Rjdefinition2}) above define $\theta(|n-\gamma|)$  closed subintervals  $Q_0,Q_1,\ldots, Q_{\theta(|n-\gamma|)-1}$  of  $[0,\pi/2]$.  We put
\begin{equation}
{\mathcal P}(\gamma)=\big\{Q_0,Q_1,Q_2,\ldots,Q_{\theta(|n-\gamma|)-1}\big\}.\label{eq:defpartition}
\end{equation}
Note that ${\mathcal P}(\gamma)$ depends only upon $|n-\gamma|$ so, strictly speaking, ${\mathcal P}(\gamma)$ depends also on $n$ as well as $\gamma$. The significance of   ${\mathcal P}(\gamma)$ lies in its relationship to the zeros of $\sin(n-\gamma)x$ in $[0,\pi/2]$. There are $\theta(|n-\gamma|)$ zeros $c_1,c_2,\ldots,c_{\theta(|n-\gamma|)}$ of  $\sin(n-\gamma)x$ in $[0,\pi/2]$ given  by
\begin{equation}
c_j=\frac{\pi j}{|n-\gamma|}, \ {\rm for}\ j=0,1,2,\ldots, \theta(|n-\gamma|)-1.\label{eq:zerossinegamma}
\end{equation}
Now, we see from (\ref{eq:Rjdefinition1}) and (\ref{eq:Rjdefinition2}) that
$c_0$ is the left endpoint of $Q_0$, when $|n-\gamma|$ is even  $c_j$ is the midpoint of $Q_j$ for $j=1,2\ldots,\theta(|n-\gamma|)-2$ and $c_{\theta(|n-\gamma|)-1} =\pi/2$ is  the right endpoint of $Q_{\theta(|n-\gamma|)-1}$ and, when $|n-\gamma|$ is odd $c_j$ is the midpoint of $Q_j$ for $j=1,2\ldots,\theta(|n-\gamma|)-1$.
\begin{lemma}\label{lemma:leavedpartition}
 Let  $\alpha,\beta, n\in {\mathbb Z}$ be such that $n\ne \alpha$ and $n\ne \beta$. Let ${\mathcal P}(\alpha,\beta)$ 
be the partition of $[0,\pi/2]$ that is the  refinement of  ${\mathcal P}(\alpha)$ and  ${\mathcal P}(\beta)$, as given by (\ref{eq:refinement}).  Then the number of intervals in ${\mathcal P}(\alpha,\beta)$ is bounded above by
\begin{equation}2\,\max{\big\{|n-\alpha|, |n-\beta|\big\}}.\label{partitionestimate}
\end{equation}
  Also, if $J\in  {\mathcal P}(\alpha,\beta)$, 
  \begin{equation}0<\lambda(J)\le\min\left\{ \frac{\pi}{|n-\alpha|}, \frac{\pi}{|n-\beta|}\right\}.\label{eq:length}
 \end{equation}
 \end{lemma}
 
 {\bf Proof.} The partition ${\mathcal P}(\alpha)$ has $\theta(|n-\alpha|)$ intervals, while ${\mathcal P}(\beta)$ has $\theta(|n-\beta|)$ intervals.   So, we see from (\ref{eq:partitionestimate1}) that  ${\mathcal P}(\alpha,\beta)$ has at most\hfill\break  $\theta(|n-\alpha|)+\theta(|n-\beta|)-1$ intervals. However, $\theta(\hbox{$r$})\le (r+2)/2$   for all $r\in {\mathbb N}$, so an upper  bound for the number of intervals in   ${\mathcal P}(\alpha,\beta)$ is 
 \[\frac{1}{2}(|n-\alpha|+|n-\beta|)+1\le\max\{|n-\alpha|, |n-\beta|\}+1\le 2\max\{|n-\alpha|, |n-\beta|\}. \]
 
 Finally, if $J\in  {\mathcal P}(\alpha,\beta)$, $J=R\cap S$ for some $R\in  {\mathcal P}(\alpha)$ and $S\in  {\mathcal P}(\beta)$.  Then, $\lambda(\hbox{$R$})\le \pi/|n-\alpha|$ and $\lambda(S)\le \pi/|n-\beta|$, and so (\ref{eq:length}) follows.
 \hfill$\square$

 Figure 1 illustrates Lemma  \ref{lemma:leavedpartition} in the case  $\alpha=1$, $\beta=-1$ and $n=9$. 

   \[\begin{tikzpicture}[scale=0.9]
\draw (0,0)--(10,0);
\draw[gray, , line width =0.12cm](0,0)--(1,0);
\draw[gray, line width =0.12cm](1,0)--(1.25,0);
\draw[gray, line width =0.12cm](1.25,0)--(3,0);
\draw[gray, line width =0.12cm](3,0)--(3.75,0);
\draw[gray, line width =0.12cm](3.75,0)--(5,0);
\draw[gray, line width =0.12cm](5,0)--(6.25,0);
\draw[gray, line width =0.12cm](6.25,0)--(7,0);
\draw[gray, line width =0.12cm](7,0)--(8.75,0);
\draw[gray, line width =0.12cm](8.75,0)--(9.0,0);
\draw[gray, line width =0.12cm](9.0,0)--(10.0,0);
      \filldraw[black] (0,0) circle (3pt);
   
     \filldraw [black](2,0) circle (3pt);
      \filldraw [black](4,0) circle (3pt);
        \filldraw [black](6,0) circle (3pt);
         \filldraw [black](8,0) circle (3pt);
         \filldraw [black](10,0) circle (3pt);
          \filldraw [black](2.5,0) circle (3pt);
           \filldraw [black](5.01,0) circle (3pt);
        
         \filldraw [black](7.5,0) circle (3pt);
         \filldraw [black](8,0) circle (3pt);
          \node[inner sep=20pt, anchor=north] at (0,0) { {$a_0=b_0$}};
  
   \node[inner sep=23pt, anchor=north] at (2.5,-0.0) {{$a_1$}};
    \node[inner sep=20pt, anchor=north] at (5.0,-0.08) {{$a_2$}};
  \node[inner sep=20pt, anchor=north] at (2.02,0) {{$b_1$}};
  \node[inner sep=20pt, anchor=north] at (4.1,0) { {$b_2$}};
   \node[inner sep=20pt, anchor=north] at (6.0,0) { {$b_3$}};
    \node[inner sep=20pt, anchor=north] at (8.0,0) { {$b_4$}};
    \node[inner sep=20pt, anchor=north] at (7.5,-0.085) { {$a_3$}};

     \node[inner sep=15pt, anchor=south] at (10,0) { {$\pi/2$}};
      \node[inner sep=20pt, anchor=north] at (10,0) { {$a_4=b_5$}};
       \node[inner sep=15pt, anchor=south] at (5,0) { {$\pi/4$}};
      \node[inner sep=15pt, anchor=south] at (0,0) { {$0$}};
   \draw (0,-0.4)--(0,0.4);
\draw (1.0,-0.4)--(1.0,0.4);
\draw (1.25,-0.4)--(1.25,0.4);
\draw (3,-0.4)--(3,0.4);
\draw (3.75,-0.4)--(3.75,0.4);
\draw (7,-0.4)--(7,0.4);
\draw (5,-0.4)--(5,0.4);
\draw (6.25,-0.4)--(6.25,0.4);
\draw (8.75,-0.4)--(8.75,0.4);
\draw (9,-0.4)--(9,0.4);
\draw (10,-0.4)--(10,0.4);
\end{tikzpicture}\]

\vskip -0.8cm
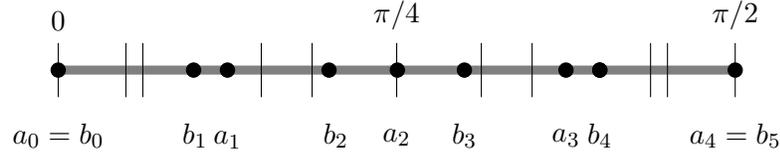
\[ \begin{minipage}{11.6cm} {\footnotesize{\bf Figure 1.} The figure illustrates the case  $n=9$, $\alpha=1$, $\beta=-1$.  We have $n-\alpha=8$, $n-\beta=10$.  The zeros  of $\sin 8x$ and $\sin 10x$  in $[0,\pi/2]$ are denoted respectively by $a_j$ and $b_j$.      The vertical lines in the figure illustrate the intervals in the refinement ${\mathcal P}(-1,1)$ of ${\mathcal P}(-1)$ and ${\mathcal P}(1)$.    Note that four of the ten intervals in ${\mathcal P}(-1,1) $ contain no zeros of   $\sin 8x\sin 10x$. }
 \end{minipage}\]
\vskip 0.3cm

Now, if $x\in {\mathbb R}$, let $d_{\mathbb Z}(x)$  denote the distance from $x$ to a nearest integer. For later use, note the fact that  $d_{\mathbb Z}(x)=|x|$ if and only if  $ -1/2\le x\le 1/2$. 
  \begin{lemma}\label{lemma:triginequality}  Let   $n, \alpha,\beta\in{\mathbb Z}$  with $n\ne \alpha$, $n\ne \beta$. Let the partitions ${\mathcal P}(\alpha)$ and ${\mathcal P}(\beta)$ be given as in (\ref{eq:defpartition}), and we write 
  $${\mathcal P}(\alpha) =\{R_0, R_1, \ldots, R_{\theta(|n-\alpha|)-1}\} \ \hbox{and}\  {\mathcal P}(\beta) =\{S_0, S_1, \ldots, S_{\theta(|n-\beta|)-1}\}.$$ 
Let  $R_j\cap S_k$ be an element of the refinement  ${\mathcal P}(\alpha,\beta)$ of  ${\mathcal P}(\alpha)$ and ${\mathcal P}(\beta)$.
    Then, if $x\in R_j\cap S_k$,  we have
 \begin{align}&\sin^2 ((n-\alpha)x)\,\sin^{2} ((n-\beta)x)\nonumber\\
 &\hskip 0cm \ge \frac{2^{4}(n-\alpha)^{2}\, (n-\beta|^{2}}{\pi^{4}}\left(x-\frac{j\pi}{|n-\alpha|}\right)^{2}\left(x-\frac{k\pi}{|n-\beta|}\right)^{2}.\label{eq:y4}
 \end{align}
 \end{lemma}
 
 {\bf Proof.} 
 We will use the fact that for all $x\in {\mathbb R}$, $|\sin \pi x|\ge 2 d_{\mathbb Z}(x)$ \cite[page 89, for example]{nillsen2}.  
    Given that $x\in R_j$, we see from (\ref{eq:Rjdefinition1}) and (\ref{eq:Rjdefinition2})  that 
 \begin{equation}
 \left|(n-\alpha)\left(\frac{x}{\pi}-\frac{j}{|n-\alpha|}\right)\right|\le\frac{1}{2}.\label{eq:y3}
 \end{equation}
 Using   (\ref{eq:y3}) we now have, for all $x$ in $R_j$,
 \begin{align} |\sin((n-\alpha)x)|&= |\sin(|n-\alpha|x-j\pi)|\nonumber\\
 &=\left|\sin\left(\pi(n-\alpha)\left(\frac{x}{\pi}-\frac{j}{|n-\alpha|}\right)\right)\right|\nonumber\\
 &\ge 2 d_{\mathbb Z}\left((n-\alpha)\left(\frac{x}{\pi}-\frac{j}{|n-\alpha|}\right)\right)\nonumber\\
 &= 2|n-\alpha|\left|\frac{x}{\pi}-\frac{j}{|n-\alpha|}\right|\nonumber\\
 &= \frac{2}{\pi}|n-\alpha|\left|x-\frac{j\pi}{|n-\alpha|}\right|.\label{eq:y5}
 \end{align}
 Using a corresponding argument, we see also that for all $x\in S_k$,  
  \begin{equation} |\sin ((n-\beta)x)|\ge  \frac{2}{\pi}|n-\beta|\left|x-\frac{k\pi}{|n-\beta|}\right|.\label{eq:y6}
 \end{equation}
 
Conclusion (\ref{eq:y4})  now follows from (\ref{eq:y5}) and (\ref{eq:y6}). \hfill
 $\square$

\section{Integral estimates  in $\boldsymbol{{\mathbb R}^m}$}
\setcounter{equation}{0}
In this section we develop  estimates for some integrals in ${\mathbb R}^m$, and an inequality between quadratics, with a view to proving Lemma \ref{lemma:cosine estimate}. 
\begin{lemma}\label{lemma:A} Let $s,m\in {\mathbb N}$ with $m\ge 4s+1$.  Then, there is a number $M>0$, depending upon $s$ and $m$ only,  such that  for all $b_1,b_2,\ldots,b_m>0$ and   for all   $(a_1,a_2,\ldots,a_m)  \in \prod_{t=1}^m[-b_t,b_t]$, 
\[\int_{{\textstyle\prod_{t=1}^m [-b_t,b_t]}}\ \frac {du_1du_2\ldots du_m}{\displaystyle\sum_{t=1}^m(u_t^2-a_t^2)^{2s}}\le M\Bigl(\max\Big\{b_1,b_2,\ldots,b_m\Big\}\Big)^{m-4s}.\]
\end{lemma}
PROOF. Clearly, we may assume that $0\le a_t\le b_t$ for all $t=1,2,\ldots,m$. Now, we have
\begin{align}
&\int_{\textstyle\prod_{t=1}^m[-b_t,b_t]}\ \frac {du_1du_2\ldots du_m}{\displaystyle\sum_{t=1}^m(u_t^2-a_t^2)^{2s}}\nonumber\\
&=2^m\int_{\textstyle\prod_{t=1}^m [0,b_t]}\ \frac {du_1du_2\ldots du_m}{\displaystyle\sum_{t=1}^m(u_t^2-a_t^2)^{2s}}\nonumber\\
&
 =2^m\int_{\textstyle\prod_{t=1}^m [0,b_t]}\ \frac {du_1du_2\ldots du_m}{\displaystyle\sum_{t=1}^m(u_t-a_t)^{2s}(u_t+a_t)^{2s}}\nonumber\\
&=2^m\int_{\textstyle\prod_{j=1}^m[-a_t,b_t-a_t]}\ \frac {dv_1dv_2\ldots dv_m}{\displaystyle\sum_{t=1}^mv_t^{2s}(v_t+2a_t)^{2s}},\label{eq:z1}
\end{align}
on putting $v_t=u_t-a_t$.
Now observe that if $v_t\ge 0$ then $v_t+2a_t\ge 0$, and that if $-a_t\le v_t\le 0$ then $v_t+2a_t\ge a_t\ge |v_t|$. Also, there is $C_m>0$ such that for all  $(v_1,v_2,\ldots,v_m)\in {\mathbb R}^m$, 
\begin{equation}\sum_{j=1}^mv_j^{4s}\ge C_m\left(\sum_{j=1}^mv_j^2\right)^{2s}.\label{eq:constantm}\end{equation}

If $r>0$, we denote the closed sphere $\{x:x\in {\mathbb R}^m\ {\rm and}\ |x|\le r\}$ by $S(0,r)$. Also, we put $b=\max\{b_1,b_2,\ldots,b_m\}$. Using  the preceding observations and  (\ref{eq:z1}) we have  
\begin{align*}  &\int_{\textstyle\prod_{t=1}^m [-b_t,b_t]}\ \frac {du_1du_2\ldots du_m}{\displaystyle\sum_{t=1}^m(u_t^2-a_t^2)^{2s}}\hskip2cm
\\
&\le2^m\int_{\textstyle\prod_{j=1}^m[-a_t,b_t-a_t]}\ \frac {dv_1dv_2\ldots dv_m}{\displaystyle\sum_{t=1}^mv_t^{4s}}\hskip2cm\\
\end{align*}
\begin{align*}
&\le\frac{2^m}{C_m}\int_{\textstyle\prod_{j=1}^m[-a_t,b_t-a_t]}\ \frac {dv_1dv_2\ldots dv_m}{\left(\displaystyle\sum_{t=1}^mv_t^2\right)^{2s}}, \hbox{using (\ref{eq:constantm})},\\
&\le\frac{2^m}{C_m} \int_{\textstyle\prod_{j=1}^m[-b_t,b_t]}\ \frac {dv_1dv_2\ldots dv_m}{\left(\displaystyle\sum_{t=1}^mv_t^2\right)^{2s}},\\    &\hskip 4.8cm{\rm  as}\ [-a_t,b_t-a_t]\subseteq [-b,b]^m,\\
\end{align*}
\begin{align*}
&\le\frac{2^m}{C_m} \int_{S(0,b{\sqrt m})}\ \frac {dv_1dv_2\ldots dv_m}{\left(\displaystyle\sum_{t=1}^mv_t^2\right)^{2s}},\  {\rm as}\ [-b,b]^m\subseteq S(0,b{\sqrt m}),\\
&=\frac{2^{m+1}\pi^{m/2}}{C_m\Gamma(m/2)}\int_0^{b{\sqrt m}}r^{m-4s-1}\,dr,\ \hbox{by\ \cite[pages 394-395]{stromberg1},}\\ 
&=\frac{2^{m+1} \pi^{m/2}m^{(m-4s)/2} b^{m-4s}}{C_m\Gamma(m/2)(m-4s)}.
\end{align*}
\noindent So, the result holds if  
$ M= C_m^{-1}2^{m+1}\pi^{m/2}\ m^{(m-4s)/2}(m-4s)^{-1}\Gamma(m/2)^{-1}.$    
${\ }$ \hfill$\square$ 
 
\begin{lemma}\label{lemma:C}  Let $s,m\in {\mathbb N}$ with $m\ge 4s+1$, and   let numbers  $b_{t,k}, c_t$ and $d_t$ be given for all $t,k$ with $t=1,2\ldots,m$ and $k=1,2$.      We assume that  
\[0\le b_{t,1}\le  c_{t}\le d_t\le b_{t,2}\]
for all $t=1,2,\ldots, m$. Then, there is a number $M>0$, depending upon $s$ and $m$ only and independent of $b_{t,k}, c_t$ and $d_t$,  such that  
\begin{align*}&\int_{\textstyle{\prod_{t=1}^m [b_{t,1}, b_{t,2}] }}\,\frac{du_1du_2\ldots du_m}{\displaystyle \sum_{t=1}^{m}\left(u_t-c_t\right) ^{2s}\left(u_t-d_t \right)^{2s}}\\
&\hskip 2cm\le M\,\Big(\max\Big\{b_{1,2}-b_{1,1},b_{2,2}-b_{2,1},\ldots,b_{m,2}-b_{m,1}\Big\} \Big)^{m-4s}. 
\end{align*}
\end{lemma}
  
{\bf Proof.}
Put, for $t=1,2,\ldots,m$,
\[\eta_t=\frac{c_t+d_t}{2}, \gamma_t=\frac{d_t-c_t}{2}, v_t=u_t-\eta_t.\]
Note that $\eta_t\in [b_{t,1},b_{t,2}]$ and $0\le \gamma_t\le (b_{t,2}-b_{t,1})/2$.
Using Lemma \ref{lemma:A} and putting $v_t=u_t-\eta_t$  in the following we have
\begin{align*}
&\int_{\textstyle\prod_{t=1}^m [b_{t,1}, b_{t,2}]}\ \frac{du_1du_2\ldots du_m}{\displaystyle \sum_{t=1}^m\left(u_t-c_t\right) ^{2s}\left(u_t-d_t \right)^{2s}}\\
&=\int_{\textstyle\prod_{t=1}^m\ [b_{t,1}-\eta_t,b_{t,2}-\eta_t]}\ \frac{dv_1dv_2\ldots dv_m}{\displaystyle \sum_{t=1}^m\left(v_t+\gamma_t\right) ^{2s}\left(v_t-\gamma_t \right)^{2s}} \\
\end{align*}
\begin{align*}
&=\int_{\textstyle\prod_{t=1}^m\ [-(\eta_t-b_{t,1}),b_{t,2}-\eta_t]}\ \frac{dv_1dv_2\ldots dv_m}{\displaystyle \sum_{t=1}^m\left(v_t^2-\gamma_t^2\right) ^{2s} }\\
&\le\int_{\textstyle\prod_{t=1}^m\ [-(b_{t,2}-b_{t,1}),b_{t,2}-b_{t,1}]}\ \frac{dv_1dv_2\ldots dv_m}{\displaystyle \sum_{t=1}^m\left(v_t^2-\gamma_t^2\right) ^{2s} },  
\end{align*}
as  $[-(\eta_t-b_{t,1}),b_{t,2}-\eta_t]\subseteq   [-(b_{t,2}-b_{t,1}),b_{t,2}-b_{t,1}]$.

We now see that Lemma \ref{lemma:A} applies because $0\le\gamma_t\le b_{t,2}-b_{t,1}$, and so there is a constant $M>0$, depending only upon $s$ and  $m$, such that
\begin{align*} \ \ \ \ \ \ \ \  \ \ \ &\int_{\textstyle\prod_{t=1}^m [b_{t,1}, b_{t,2}]}\ \frac{du_1du_2\ldots du_m}{\displaystyle \sum_{t=1}^m\left(u_t-c_t\right) ^{2s}\left(u_t-d_t \right)^{2s}}\\
  &\ \  \ \ \ \le M
\Bigl(\max\Big\{b_{1,2}-b_{1,1},b_{2,2}-b_{2,1},\ldots,b_{m,2}-b_{m,1}\Big\}\Big)^{m-4s} .\hskip 2.6cm\square
\end{align*}

\begin{lemma}\label{lemma:quadratics}Let $c\le a<b\le d$. 
Let $f$, $g$ be the quadratic functions given by $f(x)=(x-c)(d-x)$, 
$g(x)=(x-a)(b-x)$.
Then, $f(x)\ge g(x)\ge 0$ for all $a\le x\le b$.
\end{lemma}

{\bf Proof.} We have
$(f-g)(x)=x(c+d-a-b)+ab-cd.$
So,
\[(f-g)(a) =(a-c)(d-a)\ge 0\ {\rm and}\ (f-g)(b)= (d-b)(b-c)\ge 0.
\]
As $f-g$ is linear and non-negative at $a$ and $b$, we deduce that $f(x)\ge g(x)$ for all $x\in [a,b]$.
${\ }$\hfill$\square$

\section  {Completion of the proof of Theorem \ref{theorem:main2} }
 
Let $s\in {\mathbb N}$ and $n,\alpha,\beta\in {\mathbb Z}$  with $n\ne \alpha$ and $n\ne \beta$. Let  $a_0,a_1,\ldots,a_{\theta(|n-\alpha|)-1}$ and $b_0, b_1,\ldots, b_{\theta(|n-\beta|)-1}$  respectively be  the zeros of $\sin(n-\alpha)x$ and $\sin (n-\beta)x$ in $[0,\pi/2]$, as given correspondingly for $\sin(n-\gamma)x$ in (\ref{eq:zerossinegamma}). Let ${\mathcal P}(\alpha) =\{R_0, R_1, \ldots, R_{\theta(|n-\alpha|)-1}\}$ and ${\mathcal P}(\beta) =\{S_0, S_1, \ldots, S_{\theta(|n-\beta|)-1}\}$ be the partitions as given by (\ref{eq:defpartition}), and  recall that ${\mathcal A}$ is the set of all $(j,k)$ such that $\lambda(R_j\cap S_k)>0$   and that the partition ${\mathcal P}(\alpha,\beta)$ of $[0,\pi/2]$ is the set $\{R_j\cap S_k: (j,k)\in {\mathcal A}\}$, as in Lemma \ref{lemma:leavedpartition}. We see from (\ref{eq:int1})  of Lemma \ref{lemma:integralsequal} and from  (\ref{eq:y4}) of Lemma \ref{lemma:triginequality} that  for any $m,s\in {\mathbb N}$, 
\begin{align}
&\int_{{[0,2\pi)}^m}
\frac{dx_1dx_2\cdots dx_m}
{\displaystyle\sum_{j=1}^{m}\left|\cos\left(\left(\frac{(\alpha-\beta)}{2}\right)x_j\right)-\cos \left(\left(n-\frac{(\alpha+\beta)}{2}\right)x_j\right)\right|^{2s}}
 \nonumber\\
&\hskip 0.5cm=2^{2m-2s}\int_{{[0,\pi/2)}^m}\frac{dx_1dx_2\cdots dx_m}{\displaystyle\sum_{j=1}^m\sin^{2s}((n-\alpha)  x_j)\,
\sin^{2s}((n-\beta) x_j)} \nonumber\\
&\le\frac{2^{2m-6s}\pi^{4s}}{(n-\alpha)^{2s}(n-\beta)^{2s}}\sum_{(j_1,k_1),\ldots, (j_m,k_m)\in {\mathcal A}}J\bigl((j_1,k_1),\ldots, (j_m,k_m)\bigr),  \label{eq:sumintegrals}
 \end{align}
 where
 \begin{equation}J\bigl((j_1,k_1),\ldots, (j_m,k_m)\bigr)=\int_{\textstyle\prod_{t=1}^mR_{j_{t}}\cap S_{k_{t}}  }\frac{dx_1dx_2\cdots dx_m}{\displaystyle\sum_{t=1}^{m}\left(x_t-a_{j_t}\right)^{2s}\left(x_t-b_{k_t}\right)^{2s}}.\label{eq:Jintegrals}
 \end{equation}
 Note that by  Lemma \ref{lemma:leavedpartition}, the number of terms  in the sum in (\ref{eq:sumintegrals}) is bounded by 
 \begin{equation}2^m\max\{|n-\alpha|^m, |n-\beta|^m\}.\label{eq:sumbound}
 \end{equation}

Now,  let $(j_1,k_1),\ldots, (j_m,k_m)\in {\mathcal A}$, and let $t\in \{1,2,\ldots,m\}$. We consider the following possibilites (i), (ii) and (iii). 

(i) If $a_{j_t}\in R_{j_t}\cap S_{k_t}$ and $b_{k_t}\in R_{j_t}\cap S_{k_t}$  we put $a_{j_t}^{\prime}=a_{j_t}$ and  $b_{k_t}^{\prime}=b_{k_t}$. Note that if $a_{j_t}=b_{k_t}$, then  both $a_{j_t}$ and $b_{k_t}$  belong to $R_{j_t}\cap S_{k_t}$.

(ii) If   $a_{j_t}\in R_{j_t}\cap S_{k_t}$  but $b_{k_t}\notin R_{j_t}\cap S_{k_t}$, we put $a_{j_t}^{\prime}=a_{j_t}$  and   take $b_{k_t}^{\prime}$ to be the endpoint   of $R_{j_t}\cap S_{k_t}$ that is closest to $b_{k_t}$.  If   $a_{j_t}\notin R_{j_t}\cap S_{k_t}$  but $b_{k_t}\in R_{j_t}\cap S_{k_t}$, we   take $a_{j_t}^{\prime}$ to be the endpoint   of $R_{j_t}\cap S_{k_t}$ that is closest to $a_{j_t}$, and we put $b_{k_t}^{\prime}=b_{k_t}$.

(iii) If $a_{j_t}\notin R_{j_t}\cap S_{k_t}$ and $b_{k_t}\notin R_{j_t}\cap S_{k_t}$, it must happen that $a_{j_t}$ lies to the left of $R_{j_t}\cap S_{k_t}$ and $b_{k_t}$ to the right, or vice versa.  In either case, we put $a_{j_t}^{\prime}$ to be one endpoint of   $R_{j_t}\cap S_{k_t}$ and $b_{k_t}^{\prime}$ to be the other.

Now it is clear that in the cases (i) and (ii) above, 
  for all $x_t\in R_{j_t}\cap S_{k_t}$ we have
\begin{equation}|x_t-a_{j_t}|^{2s}|x_t-b_{k_t}|^{2s}\ge|x_t-a_{j_t}^{\prime}|^{2s}|x_t-b_{k_t}^{\prime}|^{2s}.\label{eq:doubleprime}\ \end{equation}
In case (iii) above, we see from the simple result on quadratics in Lemma \ref{lemma:quadratics} that (\ref{eq:doubleprime}) holds for all $x_t\in R_{j_t}\cap S_{k_t}$. All possibilities are exhausted by (i), (ii) and (iii), so  that for all $j_t, k_t$  we see that $a_{j_t}^{\prime}, b_{k_t}^{\prime}\in R_{j_t}\cap S_{k_t}$ and that (\ref{eq:doubleprime}) holds for all $x_t\in R_{j_t}\cap S_{k_t}$.

Now assume that $m\in {\mathbb N}$ with $m\ge 4s+1$. We  have from (\ref{eq:Jintegrals}) and (\ref{eq:doubleprime}) that there is $M>0$, depending on $s$ and $m$ only, such that
\begin{align}
&J\bigl((j_1,k_1),\ldots, (j_m,k_m)\bigr)\nonumber\\
&=\int_{\textstyle\prod_{t=1}^mR_{j_{t}}\cap S_{k_{t}}  }\frac{dx_1dx_2\cdots dx_m}{\displaystyle\sum_{t=1}^{m}\left|x_t-a_{j_t}\right|^{2s}\left|x_t-b_{k_t}\right|^{2s}}\nonumber\\
&\le\int_{\textstyle\prod_{t=1}^mR_{j_{t}}\cap S_{k_{t}}  }\frac{dx_1dx_2\cdots dx_m}{\displaystyle\sum_{t=1}^{m}|x_t-a_{j_t^{
\prime}}|^{2s}|x_t-b_{k_t^{\prime}}|^{2s}}, \ {\rm  by}\  (\ref{eq:doubleprime}),\nonumber\\
&\le M\left(\max\big\{\lambda(R_{j_1}\cap S_{k_1}),\lambda(R_{j_2}\cap S_{k_2}),\ldots, \lambda(R_{j_m}\cap S_{k_m}) \big\}\right)^{m-4s},\nonumber\\
&\hskip 8.5cm {\rm by\  Lemma}\  \ref{lemma:C},\nonumber\\ 
&\le\pi^{m-4s}M\min\left\{\frac{1}{|n-\alpha|^{m-4s}},\frac{1}{|n-\beta|^{m-4s}}\right\}, \, {\rm using}\  (\ref{eq:length}),\nonumber\\
&=\frac{\pi^{m-4s}M}{\max\{|n-\alpha|^{m-4s}, |n-\beta|^{m-4s}\}}.\label{eq:important}
\end{align}
Now, using (\ref{eq:sumbound}),  we  see from (\ref{eq:sumintegrals}) and (\ref{eq:important}) that
\begin{align}
&\int_{{[0,2\pi)}^m}\frac{dx_1dx_2\cdots dx_m}{\displaystyle\sum_{j=1}^{m}\left|\cos\left(\left(\frac{(\alpha-\beta)}{2}\right)x_j\right)-\cos \left(\left(n-\frac{(\alpha+\beta)}{2}\right)x_j\right)\right|^{2s}}\nonumber\\
&\le\frac{2^{3m-6s}\pi^{m}M}{(n-\alpha)^{2s}(n-\beta)^{2s}}\cdot 
 \frac{\max\{|n-\alpha|^m,|n-\beta|^m\} }{\max\{|n-\alpha|^{m-4s}, |n-\beta|^{m-4s}\}}\nonumber
 \end{align}
\begin{align}
 &= 2^{3m-6s}\pi^{m}M \cdot \max\left\{\frac{(n-\alpha)^{2s}}{(n-\beta)^{2s}}, \frac{(n-\beta)^{2s}}{(n-\alpha)^{2s}} \right\}
\nonumber \\
 &\le  2^{3m-6s}\pi^{m}MK, \label{eq:conclusion} 
\end{align}
where $K>0$ is a suitable constant chosen to be  independent of $n$.  
   Lemma \ref{lemma:cosine estimate} is immediate from  (\ref{eq:conclusion}) upon taking $m$ to be $4s+1$  and, as discussed in Section 2,  Theorem \ref{theorem:main2} is now established.\hfill$\square$
 
\section{A sharpness result}

If was shown in Theorem \ref{theorem:main2} that if $f\in L^2([0,2\pi])$ is such that ${\widehat f}(\alpha)={\widehat f} (\beta)=0$,  then for almost all $(b_1,b_2, \ldots,b_{4s+1})\in [0,2\pi]^{4s+1}$, $f$ can be written in the form (\ref{eq:differencesumb}) and consequently in the form (\ref{eq:differencesum}). However, in this section we show that if  $m\in {\mathbb N}$ and  $b_1,b_2,\ldots, b_m\in [0,2\pi]^{4s+1}$ are given, there are many functions with ${\widehat f}(\alpha)={\widehat f} (\beta)=0$ that cannot be written in the form (\ref{eq:differencesum}). Thus, no single choice of $b_1,b_2,\ldots, b_m$ suffices to ensure that $(\ref{eq:differencesum})$ is possible for all $f\in L^2([0,2\pi])$  such that ${\widehat f}(\alpha)={\widehat f} (\beta)=0$.  The methods extend the techniques in \cite[pages 420-421]{meisters1}.

\begin{lemma} \label{lemma:diophantine} Let $c_1,c_2,\ldots,c_m\in {\mathbb R}$. Then, there are infinitely many $q\in {\mathbb N}$ such that $d_{\mathbb Z}(qc_j)<1/q^{1/m}$ for all $j=1,2,\ldots,m$.
\end{lemma}
{\bf Proof.} See \cite[Theorem 4.6]{niven1} or \cite[page 27]{schmidt1}, for example.\hfill$\square$
\begin{theorem} 
 Let  $m, s\in {\mathbb N}$ and let  $\alpha,\beta\in {\mathbb Z}$ be given. Also,  let $c_1,c_2,\ldots,c_m\in [0,2\pi]$ be given. Then, there is a vector subspace $V$ of $L^2([0,2\pi])$ such that $V$  has algebraic dimension equal to that of the continuum but,  for any $f\in V$ with $f\ne 0$,  there is no choice  of $f_1,f_2,\ldots,f_m\in L^2([0,2\pi])$ such that $f$  is equal to   
 \[ \sum_{j=1}^m\ \left[
 \left( e^{ic_j\left(\frac{\alpha-\beta}{2}\right)}+e^{-ic_j\left(\frac{\alpha-\beta}{2}\right)} \right) 
 \delta_0-\left( e^{ic_j\left(\frac{\alpha+\beta}{2}\right)}\delta_{c_j}+e^{-ic_j\left(\frac{\alpha+\beta}{2}\right)} \delta_{-c_j} \right)
 \right]^s\ast f_j .
  \]  
\end{theorem}

{\bf Proof.} Let $f_1,f_2,\ldots,f_m\in L^2([0,2\pi])$, for $b\in [0,2\pi]$ let $\lambda_b$ be given by   (\ref{eq:lambdasubb}), and let $f$ be given by  
\begin{equation}f=\sum_{j=1}^m\lambda_{c_j}^s\ast f_j.\label{eq:formoff}
\end{equation}
Then, using (\ref{eq:Fouriertransform}), for all $n\in {\mathbb Z}$ we have
\begin{align*}
{\widehat f}(n)&=\sum_{j=1}^m{\widehat {\lambda}_{c_j}} (n)^s {\widehat {f_j}}(n)\\
&=\sum_{j=1}^m \left(
\cos\left(\left(\frac{\alpha-\beta}{2}\right)c_j\right) -\cos\left(\left(n-\frac{\alpha+\beta}{2}\right)c_j\right)\right)^s {\widehat {f_j}}(n)\\
&=2^s\sum_{j=1}^m  \sin^s\left(\frac{(n-\alpha)c_j}{2}\right)\sin^s\left(\frac{(n-\beta)c_j}{2}\right)\,{\widehat {f_j}}(n).
\end{align*}
Thus, as $|\sin x|\le \pi d_{\mathbb Z}(x/\pi)$, we see that
\begin{equation} |{\widehat f}(n+\alpha)|\le 2^s\sum_{j=1}^m\left|\sin \left(\frac{nc_j}{2}\right)\right|^s\, |{\widehat {f_j}}(n+\alpha)|\le \pi^s 2^s\sum_{j=1}^md_{\mathbb Z}\left( \frac{nc_j}{2\pi}\right)^s\, |{\widehat {f_j}}(n+\alpha)|.\label{eq:diophantineestimate1}
\end{equation}
 
 Now, by Lemma \ref{lemma:diophantine}, for each $\ell\in {\mathbb N}$ and $k\in \{1,2 \ldots,2^{\ell}\}$ there is $q_{k,\ell}\in {\mathbb Z}$ such that
\begin{equation}
d_{\mathbb Z}\left(q_{k,\ell} \left(\frac{c_j}{2\pi}\right)\right)<\frac{1}{\ell^{1/m}}, \ {\rm for \ all}\ j=1,2,\ldots,m.\label{eq:diophantineestimate2} 
\end{equation}
We see also from Lemma \ref{lemma:diophantine} that the integers $q_{k,\ell}$ can be chosen so that they are all distinct, over all $\ell\in {\mathbb N}$ and $k\in\{1,2,\ldots,2^{\ell}\}$.
Now, let $\Phi$ be the family of all functions $\phi:{\mathbb N}\longmapsto {\mathbb N}$ such that $\phi(1)\in \{1,2\}$ and $\phi(\ell+1)\in \{2\phi(\ell)-1, 2\phi(\ell)\}$ for all $\ell\in {\mathbb N}$.  Note that  if $\phi,\phi^{\prime}\in \Phi$ and $\phi(j)\ne\phi^{\prime}(j)$, then 
\begin{equation}\phi(n)\ne\phi^{\prime}(n),\ {\rm for\  all}\  n>j.\label{eq: notequal}
\end{equation}
 We see that $\phi(n)\in \{1,2,\ldots,2^n\}$ for all $n\in {\mathbb N}$, and that $\Phi$ has the cardinality of the continuum. Now, if $\phi\in \Phi$, define a function $f_{\phi}$ as follows.  If $n\in {\mathbb Z}$ and $n=q_{\phi(\ell),\ell}+\alpha$ for some $\ell\in {\mathbb N}$, then $\ell$ is unique and we put  
\[{\widehat {f_{\phi}}}(n)=\frac{1}{\ell^{1/2+s/m}}.\]
If $n\notin\{q_{\phi(\ell),\ell}+\alpha:\ell\in {\mathbb N}\}$, we put
\[{\widehat {f_{\phi}}}(n)=0.\]
Then, because  
\[\sum_{n=-\infty}^{\infty}|{\widehat {f_{\phi}}}(n)|^2=\sum_{\ell=1}^{\infty}\frac{1}{\ell^{1+2s/m}}<\infty,\]
we see that $f_{\phi}\in L^2([0,2\pi])$.

Now, assume that $\phi\in\Phi$   and that $f_{\phi}$ can be put in the form (\ref{eq:formoff}).
By (\ref{eq:diophantineestimate1}) and (\ref{eq:diophantineestimate2}) and taking $n=q_{\phi(\ell),\ell}$ where   $\ell\in {\mathbb N}$ we have
\[|{\widehat f}_{\phi}(q_{\phi(\ell),\ell}+\alpha)|\le \frac{\pi^s2^s}{\ell^{s/m}}\left(\sum_{j=1}^m  |{\widehat {f_j}}(q_{\phi(\ell),\ell}+\alpha)|\right)\le\frac{C\pi^s2^s}{\ell^{s/m}}\left(\sum_{j=1}^m  |{\widehat {f_j}}(q_{\phi(\ell),\ell}+\alpha)|^2\right)^{1/2},\]
for some $C>0$ that depends upon $m$ only, and we see that
\begin{equation}\sum_{\ell=1}^{\infty}\ell^{2s/m}|{\widehat f}_{\phi}(q_{\phi(\ell),\ell}+\alpha )|^2<\infty.\label{eq:contra1}
\end{equation}
However,
\begin{equation}\sum_{\ell=1}^{\infty}\ell^{2s/m}|{\widehat f}_{\phi}(q_{\phi(\ell),\ell}+\alpha)|^2=\sum_{\ell=1}^{\infty}\frac{\ell^{2s/m}}{\ell^{1+2s/m}}=\sum_{\ell=1}^{\infty}\frac{1}{\ell}=\infty.\label{eq:contra2}
\end{equation}
The contradiction between (\ref{eq:contra1}) and  (\ref{eq:contra2}) shows that $f_{\phi}$ cannot be put in the form (\ref{eq:formoff}). 

Now, let $\phi_1,\phi_2,\ldots,\phi_r\in \Phi$ be distinct, and consider a linear combination $h=\sum_{j=1}^rd_jf_{\phi_j}$ where, say, $d_1\ne 0$.  It follows from the observation (\ref{eq: notequal}) above that there is $k_0\in {\mathbb N}$ such that for all $k>k_0$ with ${\widehat f}_{\phi_1}(k)\ne 0$, we have ${\widehat f}_{\phi_j}(k)=0$ for all $j\in \{2,3,\ldots,r\}$. Thus,  for all $k>k_0$  with ${\widehat f}_{\phi_1}(k)\ne 0$ , we see that ${\widehat h}(k)=d_1{\widehat f}_{\phi_1}(k)$. Then, if we apply (\ref{eq:contra1}) and (\ref{eq:contra2})  with $f_{\phi_1}$ in place of $f_{\phi}$, we see that the contradiction between (\ref{eq:contra1}) and (\ref{eq:contra2}) applies to $h$, and we deduce that $h$ cannot be written in the form (\ref{eq:formoff}). 
We now see that if  $V$ is the subspace of $L^2([0,2\pi])$ finitely spanned by $\{f_{\phi}:\phi\in \Phi\}$, then $V$ has the required properties.  \hfill$\square$

\section{Results for compact, connected abelian groups}
\setcounter{equation}{0}
In this section we look at some applications of the earlier results  to compact, connected abelian  groups and to the automatic continuity of linear forms on $L^2$ spaces on these groups.  

{\bf Definitions.}   If $z\in {\mathbb T}$ and $\nu\in {\mathbb R}$ we may write $z$ uniquely as $z=e^{it}$ where $t\in [0,2\pi)$, and we then take  $z^{\nu}$ to be $e^{it\nu}$. In particular, for $\alpha\in {\mathbb Z}$,  $z^{\alpha/2}$ is $e^{it\alpha/2}$. Let $\alpha, \beta\in {\mathbb Z}$ and let $L$ be a linear form on $L^2({\mathbb T})$.
Then $L$ is called \emph{$(\alpha,\beta)$-invariant} if, for all $b\in {\mathbb T}$ and $f\in L^2({\mathbb T})$, 
\[L\left[\left(b^{(\alpha+\beta)/2}\delta_b+b^{-(\alpha+\beta)/2}\delta_{b^{-1}}\right)\ast f\right]= \bigl(b^{(\alpha-\beta)/2}+b^{-(\alpha-\beta)/2}\bigr)\,L(f).\]
Also, $L$ is called \emph{translation invariant} if $L(\delta_b\ast f)=L(f)$ for all $b\in G$ and $f\in L^2({\mathbb T})$.

\begin{theorem}\label{theorem:circlegroupresult} Let $\alpha, \beta\in {\mathbb Z}$ and let $L$ be a linear form on $L^2({\mathbb T})$. Then the following hold.

(i) If $L$ is $(\alpha,\beta)$-invariant on $L^2({\mathbb T})$, then $L$ is continuous on $L^2({\mathbb T})$. Thus, in this case there is a function $h\in L^2({\mathbb T})$ such that $L(f)=\int_{\mathbb T}f{\overline h}\,d\mu_{\mathbb T}$ for all $f\in  L^2({\mathbb T}$.

(ii)  If $h\in L^2({\mathbb T})$, the linear form on $h\in L^2({\mathbb T})$ given by $f\longmapsto \int_{\mathbb T}f{\overline h}\,d\mu_{\mathbb T}$ is $(\alpha,\beta)$-invariant  if and only if there are $c_1,c_2\in {\mathbb C}$ such that $h(z)=c_1z^{\alpha}+c_2z^{\beta}$ for almost all $z\in {\mathbb T}$.

(iii) If $L$ is a translation invariant linear form on $L^2({\mathbb T})$, $L$ is a multiple of the Haar measure on $\mathbb T$.

\end{theorem}
{\bf Proof.} (i) It follows immediately from Theorem \ref{theorem:main2} with $s=1$ that if $L$ is $(\alpha,\beta)$-invariant,  $L$ vanishes on the space 
\[\bigl\{f:f\in L^2({\mathbb T})\ {\rm and}\ {\widehat f}(\alpha)={\widehat f}(\beta)=0\bigr\}.\]
Thus, if $L$ is $(\alpha,\beta)$-invariant  it vanishes on this    closed subspace of $L^2({\mathbb T})$, a space that has finite codimension in $L^2({\mathbb T})$. By (d) of Proposition 5.1 in \cite[page 25]{nillsen2},  $L$ is continuous. The  existence of the function $h$ in this case comes  from the fact that  the dual space of $L^2({\mathbb T})$ is identified  with $L^2({\mathbb T})$.

(ii) If $h\in L^2({\mathbb T})$, by Theorem \ref{theorem:main2} the linear form   corresponding to $h$ is   $(\alpha,\beta)$-invariant if and only if $\int_{\mathbb T}f{\overline h}\,d\mu_{\mathbb T}=0$  whenever ${\widehat f}(\alpha)={\widehat f}(\beta)=0$. This occurs  when ${\widehat h}(n)=0$ for all $n\in {\mathbb Z}$ with $n\ne \alpha, \beta$. But  that means  that the Fourier expansion of $h$ in $L^2({\mathbb T})$  is a linear combination of $z^{\alpha}$ and $z^{\beta}$.

(iii) If $L$ is translation invariant, we see that it is $(0,0)$ invariant. Then, (ii) shows that $L$ is a multiple of the Haar measure on $\mathbb T$. \hfill$\square$
 
 Note that the conclusion (iii) in  Theorem \ref{theorem:circlegroupresult} is due to Meisters and Schmidt \cite{meisters1}.  Conclusion (ii) generalises their result.

Let $G$ denote a compact, connected  abelian group with dual group ${\widehat G}$. The identity element  in   ${\widehat G}$ is denoted   by   ${\widehat e}$. The group operation in such a  group will be written multiplicatively, and the normalised Haar measure on such a group $G$ will be denoted by $\mu_G$. We denote by $M(G)$ the family of bounded complex Borel measures on $G$.  Let $m\in {\mathbb N}$ and for each  $\gamma \in {\widehat G}$ with $\gamma \ne {\widehat e}$\, let  $h_{\gamma}:G^m\longrightarrow {\mathbb T}^m$ be the function  given by
\[h_{\gamma}(g_1,g_2,\ldots,g_m)=(\gamma(g_1), \gamma(g_2),\ldots,\gamma(g_m)).\]
 The function $h_{\gamma}$ is continuous and, as $\gamma\ne {\widehat e}$, $h_{\gamma}$  maps $G^m$ onto a compact connected subgroup of  ${\mathbb T}^m$ that is strictly larger than $\{1\}^m$, so this  connected subgroup must be ${\mathbb T}^m$ itself,  as ${\mathbb T}^m$ is connected. Consequently, $h_{\gamma}$ maps $G^m$ onto ${\mathbb T}^m$. It follows that for any non-negative measurable  function on ${\mathbb T}^m$,    we have 
\begin{equation}\int_{G^m}f\circ h_{\gamma}\,d\mu_{G^m}=\int_{{\mathbb T}^m}f\,d\mu_{\mathbb T}.\label{eq:Haarinvariance}\end{equation}
This is because each side of (\ref{eq:Haarinvariance}) defines a translation invariant integral over ${\mathbb T}^m$, so the equality in (\ref{eq:Haarinvariance})  is a consequence of the uniqueness of the Haar measure on $\mathbb T$, mentioned in (iii) of Theorem \ref{theorem:circlegroupresult}.

 In the following result, note that  for any   compact connected abelian group $G$,  for every element $b\in G$ we have $b=d^2$ for some $d\in G$ \cite[vol. I, page 385]{hewitt1}.

\begin{theorem}\label{theorem:groupsresult}
Let $G$ be a compact connected abelian group with dual group ${\widehat G}$. Let ${\widehat e}$ be the identity element of  ${\widehat G}$.
Let $s\in {\mathbb N}$ and let 
$n,\alpha,\beta\in {\mathbb Z}$ with $n\notin \{\alpha,\beta\}$. Then,   for a function $f\in L^2(G)$  the following conditions (i) and (ii)  are equivalent.

(i) The Fourier transform  of $f$ vanishes at ${\widehat e}$.  That is, ${\widehat f}({\widehat e})=0$.

(ii) There are  $b_1,b_2,\ldots,b_{4s+1}\in G$ and  $d_1,d_2,\ldots,d_{4s+1}\in G$   with  $d_j^2=b_j$ for all $j\in \{1,2,\ldots,4s+1\}$, such that  there are  $f_1,f_2,\ldots,f_{4s+1}\in L^2(G)$  so that
\begin{equation}
f=\sum_{j=1}^{4s+1}\left(\delta_{d_j^{\alpha-\beta}}+\delta_{d_j^{-(\alpha-\beta)}}-\delta_{d_j^{2n-(\alpha+\beta)}}-\delta_{d_j^{-(2n-(\alpha+\beta))}}\right)^{s}\ast f_j. \label{eq:groupidentity}
\end{equation}

When  $f$ satisfies conditions (i) and (ii),  almost all $(b_1,b_2,\ldots,b_{4s+1})\in G^{4s+1}$ have the following property:     for any $d_1,d_2,\ldots,d_{4s+1}\in G$   with  $b_j=d_j^2$ for all $j\in \{1,2,\ldots,4s+1\}$, there are $f_1,f_2,\ldots,f_{4s+1}\in L^2(G)$ such that 
(\ref{eq:groupidentity}) holds.

In the case when ${\widehat f}({\widehat e})=0$ and $\alpha-\beta$ is even, for almost all $(b_1,\ldots,b_{4s+1})\in G^{4s+1}$, there are $f_1,f_2,\ldots,f_{4s+1}\in L^2(G)$ such that
\begin{equation}
f=\sum_{j=1}^{4s+1}\left(\delta_{b_j}^{(\alpha-\beta)/2}+\delta_{b_j}^{-(\alpha-\beta)/2}-\delta_{b_j}^{n-(\alpha+\beta)/2}-\delta_{b_j}^{-(n-(\alpha+\beta)/2))}\right)^{s}\ast f_j. \label{eq:groupidentity2}
\end{equation}
\end{theorem}

{\bf Proof.}  Assume that (i) holds. 
In (\ref{eq:Haarinvariance}),  given $s\in {\mathbb N}$ and   $n,\alpha,\beta\in {\mathbb Z}$ such that $n\ne\alpha$ and $n\ne \beta$, let's take $f$ to be the function on ${\mathbb T}^{4s+1}$  whose value $f(z_1,z_2,\ldots,z_{4s+1})$ at $(z_1,z_2,\ldots,z_{4s+1})$ is
\[\frac{1}{\displaystyle\sum_{j=1}^{4s+1}\Big|z_j^{(\alpha-\beta)/2}+z_j^{-(\alpha-\beta)/2}- z_j^{n-(\alpha+\beta)/2}-z_j^{-(n-(\alpha+\beta))/2}\Big|^{2s}}.\]
Then, using (\ref{eq:Haarinvariance}) with $m=4s+1$, for each $\gamma\in {\widehat G}$ with $\gamma\ne {\widehat e}$  we have
\begin{align}&\int_{G^{4s+1}}\frac{d\mu_G(b_1) d\mu_G(b_2)\cdots d\mu_G(b_{4s+1})} {\displaystyle\sum_{j=1}^{4s+1}\Big|\gamma(b_j)^{(\alpha-\beta)/2}+\gamma(b_j)^{-(\alpha-\beta)/2}- \gamma(b_j)^{n-(\alpha+\beta)/2}- \gamma(b_j)^{-(n-(\alpha+\beta)/2)}\Big|^{2s}}\nonumber\\
&=\frac{1}{
2^{6s+1}\pi^{4s+1} }\int_{[0,2\pi]^{4s+1}}\frac{dx_1dx_2\cdots dx_{4s+1}}{\displaystyle\sum_{j=1}^{4s+1}\left|\cos \left(\left(\frac{\alpha-\beta}{2}\right)x_j\right)-\cos\left( \left(n-\left(\frac{\alpha+\beta}{2}\right)\right)x_j\right)\right|^{2s} }.
\label{eq:equationX}
\end{align}

Now let $b,d\in G$ with $d^2=b$, and let $\gamma\in {\widehat G}$.  Then, $\gamma(d)^2=\gamma(b)$, so that  if we put $\gamma(b)=e^{i\theta}$ where $\theta\in [0,2\pi)$, we have   $\gamma(d)=e^{i\theta/2}$ or $\gamma(d)=e^{i(\theta/2+\pi)}$. In the former case we have 
\begin{equation}\gamma(d)^{\alpha-\beta}=e^{i\theta(\alpha-\beta)/2}=\gamma(b)^{(\alpha-\beta)/2},
\label{eq:sign1}\end{equation}
while in the latter case we have
\begin{equation}
\gamma(d)^{\alpha-\beta}=e^{i(\alpha-\beta)\pi}e^{i\theta(\alpha-\beta)/2}=(-1)^{\alpha-\beta}\gamma(b)^{(\alpha-\beta)/2}.
\label{eq:sign2}
\end{equation}
Similarly, when  $\gamma(d)=e^{i\theta/2}$ or $\gamma(d)=e^{i(\theta/2+\pi)}$ we have, respectively,
\begin{equation} \gamma(d)^{\alpha+\beta}= \gamma(b)^{(\alpha+\beta)/2}\ \hbox{or}\ \gamma(d)^{\alpha+\beta}=(-1)^{\alpha+\beta}\gamma(b)^{(\alpha+\beta)/2}.
\label{eq:sign3}
\end{equation}
Note that in  (\ref{eq:sign1}),  (\ref{eq:sign2}) and  (\ref{eq:sign3}), $\alpha-\beta$ and $\alpha+\beta$ are both even or both odd, so $(-1)^{\alpha-\beta}$ and $(-1)^{\alpha+\beta}$ are both equal to $1$ or both equal to $-1$.

Now, for $n\in {\mathbb Z}$ and $b,d\in G$ with $d^2=b$, put 
\[\lambda_{b,d,n}=\bigl(\delta_{d^{\alpha-\beta}}+\delta_{d^{-(\alpha-\beta)}}-\delta_{d^{2n-(\alpha+\beta)}}-\delta_{d^{-(2n-(\alpha+\beta))}}\bigr)^s\in M(G).\]
Then, for $\gamma\in {\widehat G}$,  
\begin{equation}{\widehat \lambda_{b,d,n}}(\gamma)=\bigl(\gamma(d)^{-(\alpha-\beta)}+\gamma(d)^{\alpha-\beta}- \gamma(d)^{-(2n-(\alpha+\beta))}- \gamma(d)^{2n-\alpha-\beta}\bigr)^s.
\label{eq:bdequation}\end{equation}
In view of (\ref{eq:sign1}),  (\ref{eq:sign2}) and  (\ref{eq:sign3}), we see that
\begin{equation} |{\widehat \lambda_{b,d,n}}(\gamma)|=\big|\gamma(b)^{(\alpha-\beta)/2}+\gamma(b)^{-(\alpha-\beta)/2}-\gamma(b)^{n-(\alpha+\beta)/2}-\gamma(b)^{-(n-(\alpha+\beta)/2)}\big|.
 \label{eq:FTequalszero}
 \end{equation}

  Now, for each $b\in G$, let $d_b\in G$ be any element such that  $d_b^2=b$.    As $n\notin\{\alpha,\beta\}$,  if  $M$ is the constant as in Lemma \ref{lemma:cosine estimate} and   we use (\ref{eq:equationX})   and  (\ref{eq:FTequalszero}), upon changing the order of summation and integration we have 
  \begin{align}
&\int_{G^{4s+1}}\left(\sum_{\gamma\in {\widehat G}, \gamma\ne {\widehat e}}\,\frac{|{\widehat f}(\gamma)|^2}{\displaystyle\sum_{j=1}^{4s+1}|{\widehat \lambda_{b_j,d_{b_j},n}}(\gamma)|^2}\right)\,\prod_{j=1}^md\mu_G(b_j)
\le \frac{M}{2^{6s+1}\pi^{4s+1}}\sum_{\gamma\in {\widehat G}}|{\widehat f}(\gamma)|^2,\label{eq:product}
 \end{align}
which is finite by Plancherel's Theorem (see \cite[vol. II, page 226]{hewitt1}). We deduce that provided ${\widehat f}({\widehat e})=0$,  for almost all $(b_1,b_2,\ldots,b_{4s+1})\in G^{4s+1}$ we have that \begin{equation}\sum_{\gamma\in {\widehat G} }\,\frac{|{\widehat f}(\gamma)|^2}{\displaystyle\sum_{j=1}^{4s+1}|{\widehat \lambda_{b_j,d_{b_j},n}}(\gamma)|^{2}}<\infty.\label{eq:AE}
\end{equation}
Then, (ii)  follows from (\ref{eq:bdequation}), (\ref{eq:AE}) and Theorem \ref{theorem:characterisation}, so (i) implies (ii).

It is clear from (\ref{eq:bdequation}) that if $f$ has the form (\ref{eq:groupidentity}), then ${\widehat f}({\widehat e})=0$. Thus, (ii) implies (i).

 Above, the   observation was made  that (\ref{eq:AE})  holds for almost all \hfill\break$(b_1,b_2,\ldots,b_{4s+1})\in G^{4s+1}$. If $(b_1,b_2,\ldots,b_{4s+1})$ is any such point, we deduce that for any  choice of   elements $d_1,d_2,\ldots,d_{4s+1}$ such that $d_j^2=b_j$ for all $j$,  and  when ${\widehat f}({\widehat e})=0$, we will have  (\ref{eq:groupidentity}) holding. 
  
 When ${\widehat f}({\widehat e})=0$ and $\alpha-\beta$ is even, the final conclusion derives from the above arguments and the fact that 
 \[d_j^{\alpha-\beta}=(d_j^2)^{(\alpha-\beta)/2}=b_j^{(\alpha-\beta)/2}\ \hbox{and}\ d_j^{\alpha+\beta}=b_j^{(\alpha+\beta)/2}.\]

\hfill$\square$

The following is a result concerning  automatic continuity on groups. It  is derived from Theorem \ref{theorem:groupsresult}, but only a special case is stated.  More general results can be derived from Theorem \ref{theorem:groupsresult}.

\begin{theorem}\label{theorem:Haar}
Let $G$ be a compact connected abelian group. Then the following conditions (i), (ii) and (iii) on a linear form $L:L^2(G) \longrightarrow {\mathbb C}$ are equivalent.

(i) $L$ is translation invariant.  That is, $L(\delta_g\ast f)=L(f)$, for all $g\in G$ and $f\in L^2(G)$.

(ii) There is $n\in {\mathbb Z}$ with $n\notin \{-1,1\}$  such that  
\[L\bigl((\delta_g+\delta_{g^{-1}})\ast f\bigr)=L\bigl((\delta_{g^n}+\delta_{g^{-n}})\ast f\bigr),\]
for all $g\in G$ and $f\in L^2(G)$.

(iii)   $L$ is a multiple of the Haar measure.

Also, the normalised Haar measure $\mu_{\mathbb G}$ on $G$ is unique.
\end{theorem}
{\bf Proof.} (i)  implies that (ii) holds for any $n\in {\mathbb N}$, so it must hold for any particular $n$. Also (iii) implies (i) because the Haar measure is translation invariant. Finally, assume that (ii) holds for some  $n\in {\mathbb Z}$ with $n\notin \{-1,1\}$. By Theorem \ref{theorem:groupsresult} with $s=1$, $\alpha=1$ and $\beta=-1$, we deduce that $L$ vanishes on the  closed subspace $\{f:f\in L^2(G)\ {\rm and}\ {\widehat f}({\widehat e})=0\}$. This latter space has codimension $1$ in $L^2(G)$ and so it follows easily that $L$ is continuous and is a multiple of the Haar measure on $G$ (see \cite[page 415]{meisters1}), and (iii) follows. Finally, as the Haar measure $\mu_G$ defines a translation invariant linear form on  $L^2(G)$, the equivalence of (i) and (iii) implies the uniqueness of the Haar measure.
\hfill$\square$

\noindent  Rodney Nillsen 

\noindent School of Mathematics and Applied Statistics
 
\noindent University of Wollongong

\noindent  New South Wales

\noindent AUSTRALIA 2522

\noindent email: nillsen@uow.edu.au

\noindent web page: http://www.uow.edu.au/$\sim$nillsen
\end{document}